\documentclass [11pt]{article}
\usepackage{latexsym}
\usepackage[dvips]{epsfig}

\newcommand{\qed}{\mbox{$\Diamond$}\vspace{\baselineskip}}

\newtheorem{theorem}{Theorem}

\newtheorem{lemma}{Lemma}
\newtheorem{proposition}{Proposition}
\newtheorem{definition}{Definition}
\newtheorem{corollary}{Corollary}
\newtheorem{example}{Example}
\newtheorem{observation}{Observation}

\newenvironment{proof}{\noindent {\bf Proof:}}{{\qed}}

\newcommand{\vanish}[1]{}
 
\begin{document}
\title{A self-dual poset on objects counted by the Catalan numbers
      \\and a type-B analogue}
 
\author{Mikl\'os B\'ona\thanks{This paper was written while the author's stay at 
the
Institute was supported by Trustee Ladislaus von Hoffmann, the Arcana
Foundation.}\\
        School of Mathematics\\
        Institute for Advanced Study\\
        Princeton, NJ 08540 \\
        bona@ias.edu\\
\\
\\
	Rodica Simion\\
	Department of Mathematics\\
	The George Washington University\\
	Washington, DC 20052\\
	simion@gwu.edu\\
	}

\maketitle
\begin{abstract}
We introduce two partially ordered sets, $P^A_n$ and $P^B_n$, of the same
cardinalities as the type-A and type-B noncrossing partition lattices.
The ground sets of $P^A_n$  and $P^B_n$ are subsets of the symmetric and the
hyperoctahedral groups, consisting of permutations which 
avoid certain patterns.  The order relation is given by
(strict) containment of the descent sets.  In each case, by means of an
explicit order-preserving bijection, we show that the poset of restricted
permutations is an extension of the refinement order on noncrossing
partitions.  Several structural properties of these permutation posets
follow, including  self-duality and the strong Sperner property.
We also discuss posets $Q^A_n$ and $Q^B_n$ similarly associated with
noncrossing partitions, defined by means of the excedence sets of suitable 
pattern-avoiding subsets of the symmetric and hyperoctahedral groups.

\


\end{abstract}

\section{Introduction}\label{sec-Intro}
There are more than 150 different objects enumerated by Catalan numbers; 
\cite{EC2} contains an extensive list of such combinatorial objects and
their properties. Two
of the most carefully studied ones are noncrossing partitions  and 
permutations avoiding a 3-letter pattern. 

A partition $\pi$ of the set $[n] \colon =\{1,2,\cdots, n\}$, having blocks
$\beta_1,\beta_2,\cdots ,\beta_k$, is called {\em noncrossing} 
 if there are no four elements $1 \leq a<b<c<d \leq n$ so 
that $a,c\in \beta_i$ and $b,d \in \beta_j$
for some distinct blocks  $\beta_i$ and $\beta_j$. 
The set of noncrossing partitions of $[n]$ constitutes a lattice under the
refinement order (where $\pi < \nu$ if each block of $\nu$ is a union of
blocks of $\pi$).  An investigation of structural and enumerative
properties of this lattice was initiated by Kreweras \cite{kreweras},
and continued by several authors, e.g., \cite{Ed1}, \cite{Ed2}, 
\cite{EdSi}, \cite{Montenegro}, \cite{NicaSpei}, \cite{SiU},
\cite{St-NCact}.  
We denote the lattice of noncrossing partitions of
$[n]$ as $NC^A_n$, since it is a subposet (indeed, a sub-meet-semilattice) 
of the intersection lattice associated with the type-A hyperplane
arrangement in ${\bf R}^n$ (which consists of the hyperplanes with equations
$x_i = x_j$, for $1 \leq i < j \leq n$).
For our purposes, recall from \cite{kreweras} 
that the poset $NC^A_n$ is ranked with rank
function ${\rm rk}(\pi) = n - {\rm bk}(\pi)$ (where ${\rm bk}(\pi)$ 
denotes the number of blocks of the partition $\pi$), rank-symmetric and
rank-unimodal with rank sizes given by the Narayana numbers 
$\left(  {1 \over n} {n \choose k} {n \choose {k+1}} \right)_{0 \leq k <n}$.
Furthermore, it is  
self-dual (see \cite{kreweras}, \cite{SiU}) and has the strong  
Sperner property (see \cite{SiU};  that is, for every $k$, the maximum 
cardinality of the union of $k$ antichains is the sum of the $k$ largest
rank-sizes).

A permutation 
$\sigma = \sigma_1 \sigma_2 \cdots \sigma_n$ of $[n]$, or, in what follows,
an $n$-permutation, is called {\em 132-avoiding}
 if there are no three positions $1 \leq a<b<c \leq n$ so that 
 $\sigma_a < \sigma_c < \sigma_b$.
Classes of restricted permutations avoiding other patterns are defined
similarly.  Such classes of permutations arise naturally in theoretical
computer science in connection with sorting problems (e.g., \cite{Knuth},
\cite{Tarjan}), as well as in the context of combinatorics related to
geometry (e.g., the theory of Kazhdan-Lusztig polynomials \cite{Brenti-KL}
and Schubert varieties \cite{Billey}).  
The investigation of classes of pattern-avoiding  
permutations from an enumerative and algorithmic point of view includes 
\cite{Barc}, \cite{Bona}, \cite{ChowWest}, \cite{DuGiWe}, 
\cite{NoonanZeil}, \cite{schmidt}, \cite{West},
to name a few.

In Section \ref{sec-A} of 
this paper we introduce the partially ordered set $P^A_n$ 
whose elements are the 132-avoiding $n$-permutations, ordered by 
$\sigma < \rho$ if ${\rm Des}(\sigma) \subset {\rm Des}(\rho)$,
where ${\rm Des}$ denotes the descent set of a permutation.
One can think of $P^A_n$ as a Boolean algebra of rank ${n-1}$ in which 
each element $S$ is replicated as many times as there are 132-avoiding
permutations with $S$ as the descent set. 
We show that this poset of restricted permutations is an extension
of the lattice of noncrossing partitions $NC^A_n$ by exhibiting a natural
order-preserving bijection from the dual order $(NC^A_n)^*$ to the poset 
$P^A_n$.  
This yields the fact  that $P^A_n$ has the
same rank-generating function as $NC^A_n$ (implicit in \cite{Si-ncstats}, 
where the joint distribution of the descent and major index statistics on 
132-avoiding permutations is shown to agree with the joint distribution of
the block and {\rm rb} statistic on noncrossing partitions).
It then follows that 
$P^A_n$ is rank-unimodal, rank-symmetric and strongly Sperner. 
We also prove that $P^A_n$ is itself a self-dual poset.

We also present type-B analogues of these results.  These constitute Section
\ref{sec-B} of the paper.
The notion of a type-B noncrossing partition of $[n]$ is that
first considered by Montenegro \cite{Montenegro}, 
systematically studied by Reiner \cite{Reiner-NCB},
	and further investigated by Hersh \cite{Hersh-NC}. 
These authors show that the type-B noncrossing partitions of $[n]$ form a
lattice, $NC^B_n$, which shares
naturally a variety of properties of $NC^A_n$.  In particular, $NC^B_n$ 
is a rank-unimodal, self-dual, strongly Sperner poset.   
We define a poset $P^B_n$ into which $NC^B_n$ can be
embedded via an order-preserving bijection, with properties analogous to
those obtained for type-A.  The parallel between the type-A and type-B cases
includes the fact that the poset $P^B_n$ is defined in terms of
pattern-avoiding elements of the hyperoctahedral group (or signed
permutations), ordered by
containment of the descent set.  The relevant pattern restriction 
is the simultaneous avoidance of the patterns $21$ and
${\overline{2}}\ {\overline{1}}$.
This class of restricted signed permutations
was considered in \cite{Si-Bstats}, where B-analogues are proposed for 
type-A results in \cite{Si-ncstats}
concerning combinatorial statistics for noncrossing partitions and 
restricted permutations.

In brief, a class of partitions and one of permutations are equinumerous, 
and further, the count of the partitions by number of blocks agrees with the
count of permutations by number of descents.   A similar situation arises
for certain type-B analogues of these objects.  Our results show that these
enumerative relations are manifestations of structural relations between
partial orders which can be defined naturally on the objects under
consideration.  We also discuss posets $Q^A_n$ and $Q^B_n$ of restricted 
permutations and signed permutations ordered by containment of their sets of
excedences.
The final section of the paper consists of remarks and problems for further
investigation.

\section{The type-A case}\label{sec-A}   

\subsection{A bijection and its properties}
It is not difficult to find a bijection from the set of noncrossing partitions
of $[n]$ onto that of 132-avoiding $n$-permutations. Here we 
exhibit and analyze the structure of such a bijection, $f$, which will serve as
the main tool in proving the results of this section. 
To avoid confusion, integers
belonging to a partition will be called {\em elements}, while integers 
belonging to a permutation will be called {\em entries}. An $n$-permutation
will always be written in the one-line notation, $p=p_1p_2\cdots p_n$, 
with $p_i = p(i)$ denoting its $i$th entry.

Let $\pi \in NC^A_n$.  We construct the 132-avoiding
permutation $p=f(\pi)$ corresponding to it as follows. 
Let $k$ be the largest element of $\pi$ which is in
the same block of $\pi$ as 1. Put the entry $n$ of $p$ in the $k$th 
position, i.e., set $p_k=n$. As $p$ is to be 132-avoiding, 
this implies that the entries larger
than $n-k$ are on the left of $n$ in $p$, 
and the entries smaller than or equal to $n-k$ are on the right of $n$.
Delete $k$ from $\pi$ and 
apply this procedure recursively, with obvious minor adjustments, 
to the restrictions of $\pi$ 
to the sets $\{ 1, \dots, k-1 \}$ and $\{ k+1, \dots , n \}$, which are also
noncrossing partitions. 
Namely, if  $j$  is the largest element in the same block as $k+1$, 
we set $p_{j} = n-k$, so that the restriction $\pi_1$ of 
$\pi$ to $\{k+1,k+2,\dots ,n\}$ yields a 132-avoiding permutation
of $\{1,2,\dots, n-k\}$ placed on the right of $n$ in $p = f(\pi)$.
Similarly, if in the restriction $\pi_2$ of $\pi$ to the set 
$\{ 1, 2, \dots, k-1 \}$ the largest element in the same block as 1 is equal
to $j$, we set $p_j = n-1$.  Thus, recursively,  $\pi_2$ 
yields a 132-avoiding permutation which we realize
on the set $\{ n-k+1, n-k+2, \dots, n-1 \}$ and we place it to the left of
$n$ in $p = f(\pi)$.
In other words, with a slight abuse of notation, $f(\pi)$
is the concatenation of $f(\pi_2)$, $n$, and $f(\pi_1)$, where $f(\pi_2)$
 permutes
the set $\{n-k+1,n-k+2,\cdots ,n-1\}$  and $f(\pi_1)$ permutes the set $[n-k]$.

To see that this is a bijection note that we can recover the maximum of 
the block containing the element  1 from the position of the entry $n$ in $p$,
and then proceed recursively. 

\begin{example} {\em If $\pi=(\{1,4,6\}, \{2,3\}, \{5\}, \{7,8\})$, then 
$f(\pi)=64573812$. }\end{example}

\begin{example} {\em If $p=(\{1,2,\cdots ,n\})$, then 
$f(p)=12\cdots n $. }\end{example}

\begin{example} {\em If $p=(\{1\},\{2\},\cdots ,\{n\})$, then 
$f(p)=n\cdots 21 $. }\end{example}

The following definition is widely used in the literature.

\begin{definition} 
Let $p=p_1p_2\cdots p_n$ be an $n$-permutation. We say that $i \in [n-1]$
is a {\em descent} of $p$ if $p_i>p_{i+1}$. The set of all descents of $p$
is called the {\em descent set} of $p$ and is denoted ${\rm Des}(p)$. 
\end{definition}

Now we are in a position to define the poset $P^A_n$ of 132-avoiding
permutations we want to study. 

\begin{definition} 
Let $p$ and $q$ be two 132-avoiding $n$-permutations. 
We say that $p<q$ in $P^A_n$ if 
${\rm Des}(p)\subset {\rm Des}(q)$. 
\end{definition}

Clearly, $P^A_n$ is a poset as inclusion is transitive. 
The Hasse diagram of $P^A_4$ is shown in Figure 1.

\begin{figure}[ht]
 \begin{center}
  \epsfig{file=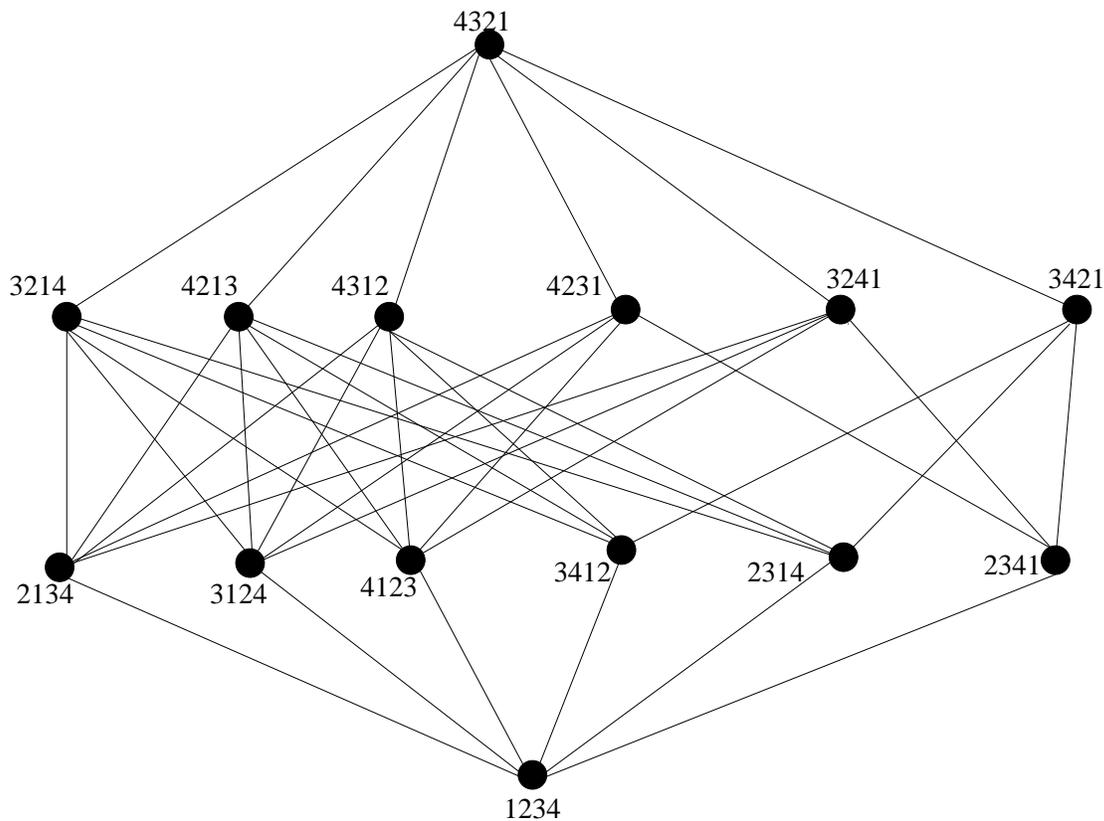}
  \caption{The Hasse diagram of $P^A_4$ }
  \label{hasse}
 \end{center}
\end{figure}

\begin{observation}\label{obs-lrminA}
In a 132-avoiding permutation, $i$
 is a descent if and only if $p_{i+1}$ is smaller
than every entry on its left.  Such an element is called a {\em left-to-right
minimum}. So $p<q$ in $P^A_n$ 
if and only if the set of positions of left-to-right minima in $p$ 
is a proper subset of the set of positions of left-to-right minima in $q$.
\end{observation}

The following proposition describes the relation between the blocks of
$\pi \in NC^A_n$ and the descent set of the 132-avoiding permutation $f(\pi)$. 

\begin{proposition} \label{prop-block} 
The bijection $f$ has the following property:  Let $i \ge 1$.  Then 
$i\in {\rm Des}(f(\pi))$ if and only if $i+1$ is the smallest element of 
its block in $\pi \in NC^A_n$.
\end{proposition}

\begin{proof} 
For $n=1$ and $n=2$ the statement is clearly true and we use induction on
$n$.  Suppose
we know the statement for all positive integers smaller than $n$. 
Then we distinguish two cases:
\begin{enumerate} \item If 1 and $n$ are in the same block of $\pi$, then the 
construction of $f(\pi)$  starts by putting the entry $n$ in the
last slot of $f(\pi)$, then deleting the element  $n$ from $\pi$. This does not
alter either the set of minimum elements of the blocks nor the set  of 
descents. Therefore, this case reduces to the general case for $n-1$, 
and is settled by the inductive hypothesis.
 \item If the largest element $k$ of the block containing 1 is smaller than
$n$, then as we have seen above,  $f(\pi)$ is the concatenation of 
$f(\pi_2), n, f(\pi_1)$, and $f(\pi_1)$ is not empty. 
Clearly, by the definition of $f(\pi)$, 
$k \in {\rm Des}(f(\pi))$ and the element 
$k+1$ is the minimum of its block.  From this and the inductive hypothesis 
applied to $f(\pi_1)$ and $f(\pi_2)$, the proof follows.
\end{enumerate}
\end{proof}

\subsection{Properties of $P^A_n$}

Proposition \ref{prop-block} implies that $P^A_n$ is isomorphic to the order 
on 
noncrossing partitions in which $\pi<\pi'$ if the set of minima of the
blocks of $\pi'$ is contained in the set of minima of the blocks of $\pi$.
This yields the first result of this section.

\begin{theorem}\label{thm-embPA}
The lattice of noncrossing partitions $NC^A_n$ is a subposet of $P^A_n$.
\end{theorem}

\begin{proof}
We show that our bijection $f$ is an order-reversing map 
$NC^A_n \to P^A_n$.  The conclusion then follows from the self-duality of
the lattice of noncrossing partitions.
Suppose $\pi < \tau$ in $NC^A_n$.  
This means $\pi$ is a finer partition than $\tau$, so every element which is
the minimum of its block in $\tau$ is also the minimum of its block in 
$\pi$.  By Proposition
\ref{prop-block} this implies ${\rm Des}(f(\tau))\subset {\rm Des}(f(\pi))$,
so $f(\pi) > f(\tau)$ in $P^A_n$.
\end{proof}

Clearly, 
$P^A_n$ is a ranked poset (with rank function 
${\rm rk}_{P^A_n}(p) = \# {\rm Des}(p)$), and we have 
${\rm rk}_{NC^A_n}(\pi) = n-1- {\rm rk}_{P^A_n}(f(\pi))$.

\begin{corollary}\label{cor-symA} 
 The poset  $P^A_n$ is rank-symmetric, rank-unimodal and strongly Sperner,
 and its rank generating function is equal to that of  $NC^A_n$.
\end{corollary}

\begin{proof}
The properties of the rank sizes of $P^A_n$ 
are immediate consequences of Proposition \ref{prop-block}
 and the corresponding properties known to hold for $NC^A_n$. 
Moreover, every antichain of $P^A_n$ is, via the bijection $f$, 
an antichain of $NC^A_n$, and the strong Sperner property of $P^A_n$ follows
from the strong Sperner property of $NC^A_n$.
\end{proof}

We now turn to showing that $P^A_n$ is self-dual, based on the next lemma.

For $S \subseteq [n-1]$, let $Perm_n(S)$ denote the
number of 132-avoiding $n$-permutations with descent set $S$. 

\begin{lemma}\label{lemma-rev-compl}
Let $S$ be any subset of $[n-1]$ and let $\alpha (S)$ denote 
its ``reverse complement,''  that is,
 $i\in \alpha(S) \Longleftrightarrow n-i \notin S$.
Then $Perm_n(S)=Perm_n(\alpha(S))$.
\end{lemma} 

\begin{proof} 
We use induction on $n$. For $n=1,2,3$ the statement is true.
Now suppose we know it for all positive integers smaller than $n$. 
Denote by $t$ the smallest element of $S$, and let $p$ be a 132-avoiding
$n$-permutation whose descent set is $S$.
\begin{enumerate} \item Suppose that $t>1$. 
Then we have $p_1<p_2<\cdots <p_t$ and,  
because $p$ avoids the pattern 132, the values of 
$p_1,p_2,\dots,p_t$ are {\em consecutive} integers. 
So, for given values of $p_1$ and $t$,
 we have only one choice for $p_2,p_3,\dots ,p_t$. This implies 
\begin{equation} \label{1a}
Perm_n(S)=Perm_{n-(t-1)}(S-(t-1)),\end{equation} 
 where $S-(t-1)$ is the set obtained from
$S$ by subtracting $t-1$ from each of its elements.

On the other hand, we have $n-t+1,n-t+2,\cdots ,n-1 \in \alpha(S)$, 
meaning that in any permutation $q$ counted by $Perm_n(\alpha(S))$ the chain
of inequalities $q_{n-t+1}>q_{n-t+2}>\cdots >q_n$ holds.  
To avoid forming a 132-pattern in $q$, we must have 
$(q_{n-t+2},\dots ,q_n)=(t-1,t-2,\dots 1)$.
Therefore, 
\begin{equation}\label{1b}
Perm_n(\alpha(S))=Perm_{n-(t-1)}(\alpha(S)|n-(t-1))
\end{equation}
where $\alpha(S)|n-(t-1)$ denotes the set obtained from  $\alpha(S)$ 
by removing its last $t-1$ elements. Clearly, 
$Perm_{n-(t-1)}(S-(t-1))=Perm_{n-(t-1)}(\alpha(S)|n-(t-1))$
by the induction hypothesis, so equations (\ref{1a}) and  (\ref{1b}) imply
 $Perm_n(S)=Perm_n(\alpha(S))$.
\item If $t=1$, but $S\neq [n-1]$,
 then let $u$ be the smallest index which is {\em not}
 in $S$. Then again, to avoid forming a 132-pattern,
 the value of $p_u$ must be the smallest positive integer $a$ which
is larger than $p_{u-1}$ and is not equal to any $p_i$ for $i\leq u-1$.
So again, we have only one choice for $p_u$.
On the other hand, the largest index in $\alpha(S)$ will be $n-u$.  Therefore,
in permutations  $q$  counted by $Perm_n(\alpha(S))$,
 we must have $q_{n-u}=1$ as Observation 1 implies that $q_{n-u}$ must be
the rightmost left-to-right minimum in such permutations,
 and that is always the entry 1.

In order to use this information to reduce our permutations in size, we define
$S'\subset [n-2]$ as follows: $i\in S'$ if and only if either $i<u$ and  then,
by the definition of $u$, 
 $i\in S$, or $i>u$ and $i+1\in S$. In other words, we decrease elements
larger than $u$ by 1; intuitively, we remove $u$ from $[n-1]$, and translate
the interval on its right one notch  to the left. If we now take $\alpha(S')$,
 that
will consist of entries $j$ so that  $j<n-u$  and
$(n-1)-(j-1)=n-j\notin S$.  So in other words, we simply 
remove $n-u$ from  $[n-1]$ (there has been nothing on the right of $n-u$ in
$\alpha(S)$ to translate). Note that
the size of $\alpha(S)$ decreases with this operation as $n-u\in \alpha(S)$.
As we have seen in the previous paragraph, we had only one choice for $p_u$
and $p_{n-u}$, so removing them this way does not change the number of
permutations with a given descent set. Thus we have  
$Perm_n(S)=Perm_{n-1}(S')$, and also $Perm_n(\alpha(S))=
Perm_{n-1}(\alpha(S'))$. By induction hypothesis, the right hand sides of
thes two equations agree, and therefore the left hand sides must agree, too.
\begin{example} {\em If $n=8$ and $S=\{1,6\}$, and so $\alpha(S)=
\{1,3,4,5,6\}$, then $u=2$, $n-u=6$, and indeed, $S'=\{1,5\}$ and $\alpha(S')=
\{1,3,4,5\}$. }
\end{example}
\item Finally, if $S=[n-1]$, then the statement is trivially true as 
$Perm_n(S)=Perm_n(\alpha(S))=1.$
\end{enumerate}
So we have seen that $Perm_n(S)=Perm_n(\alpha(S))$ in all cases.
\end{proof}

It is now easy to verify that  the reverse complementation of the 
descent set can be used to construct an anti-automorphism of $P^A_n$.

\begin{theorem} The poset $P^A_n$ is self-dual.
\end{theorem}

\begin{proof}
It is clear that, in $P^A_n$, 
permutations which have the same descent set will
cover the same elements and they will be covered by the same elements.
The permutations with a prescribed descent set $S$ 
form an orbit of $Aut(P^A_n)$ and they can be
permuted among themselves arbitrarily by elements of $Aut(P^A_n)$. 
On the other hand, Lemma \ref{lemma-rev-compl} shows that the orbits
corresponding to $S \subseteq [n-1]$ and to its reverse-complement 
$\alpha(S)$ are equinumerous.  Hence, a map $P^A_n \to P^A_n$ which
establishes a bijection between 
$\{ p \in P^A_n \ \colon \ {\rm Des}(p) = S \}$ 
and $\{ q \in P^A_n \ \colon \ {\rm Des}(q) = \alpha(S) \}$ for each $S
\subseteq [n-1]$ provides an order-reversing bijection on $P^A_n$.
\end{proof}

\subsection{A poset derived from excedences}\label{subsec-excA}

It is shown in \cite{Si-ncstats} that the joint distribution
of the excedence and Denert statistics on 321-avoiding permutations agrees
with the joint distribution of the block and ${\rm rb}$ statistics on
noncrossing partitions.   This suggests the definition of the poset 
$Q^A_n$ consisting of the 321-avoiding $n$-permutations ordered by
containment of the set of excedences, and invites the question of how 
$Q^A_n$ compares with the poset $P^A_n$.

A permutation $\sigma$ has an {\em excedence} at $i$ if $\sigma(i) > i$.
For example, the {\em excedence set} of $\sigma = 3 2 5 1 4$ is 
${\rm Exc}(\sigma) = \{ 1, 3 \}$.  Let ${\rm exc}(\sigma)$ denote the number
of excedences of $\sigma$. 
Following \cite{Si-ncstats}, 
there is a bijection $\theta$ from $NC^A_n $ to 321-avoiding $n$-permutations 
such that ${\rm exc}(\theta(\pi)) = {\rm bk}( \pi ) - 1$. 
Namely, if the set of minima of the blocks of $\pi \in NC^A_n$, 
omitting the block containing 1, is 
$\{ f_2 < \cdots < f_k \}$ and the set of maxima of the blocks, 
again, omitting the block containing 1, is
$\{ l_2 < \cdots < l_k \}$, then let $\theta(\pi)$ be the permutation 
whose value at $f_i - 1$ is $l_{i} $ for $i = 2, 3, \dots, k$, and
whose other values constitute an increasing subsequence in the remaining
positions.   For instance, if $\pi = \{ 1, 5, 7 \} \{ 2 \} \{ 3, 4 \} 
\{ 6 \} \{ 8, 10 \} \{ 9 \} \in NC^A_{10}$, then we have 
$( f_2, \dots, f_6) = ( 2, 3, 6, 8, 9)$ and 
$( l_2, \dots, l_6) = (2, 4, 6, 9, 10)$, and we obtain 
$\theta(\pi) = 2\ 4\ 1\ 3\ 6\ 5\ 9\ 10\ 7\ 8$.  

Recall from \cite{Si-ncstats} that the set of excedences of $\theta(\pi)$ is
precisely $\{ f_2 -1, f_3 -1, \dots, f_k -1 \}$.  Similarly to the case of 
descents discussed for 132-avoiding permutations, a covering relation 
$\pi < \pi'$ in $NC^A_n$ corresponds to the deletion of an excedence: 
${\rm Exc}(\theta(\pi')) = {\rm Exc}(\theta(\pi)) - \{ i \}$, for  
a suitable $i \in {\rm Exc}(\theta(\pi))$.  Hence, taking advantage
of the self-duality of $NC^A_n$, one can establish directly that 
the poset $Q^A_n$ enjoys the same properties as $P^A_n$:
There is an embedding of $NC^A_n$ into the poset $Q^A_n$ of 321-avoiding
$n$-permutations ordered by containment of the set of excedences;
the embedding is rank-preserving and $Q^A_n$ is a strongly Sperner poset.

The fact that the posets $P^A_n$ and $Q^A_n$ have strongly similar
 properties is not  accidental.

\begin{proposition}\label{prop-isoA} The posets $P^A_n$ and $Q^A_n$ are
isomorphic. 
\end{proposition}

\begin{proof} 
For each $S \subseteq [n-1]$, let $E^{321}_n(S)$ be the set of 321-avoiding 
$n$-permutations with excedence set $S\subseteq [n-1]$. 

Let also $D^{132}_n(\alpha(S))$ be the set of 132-avoiding $n$-permutations 
with descent set equal to $\alpha(S)$, the reverse-complement of $S$.
Thus, in the notation of the previous subsection, 
the cardinality of $D^{132}_n(\alpha(S))$ is $Perm_n(\alpha(S))$.

We construct a bijection $s \colon E^{321}_n(S) \to D^{132}_n(\alpha(S))$
(illustrated by example \ref{eg-sbij}).
If $p\in E^{321}_n(S)$, then, as seen earlier in the definition of $\theta$, 
the entries $p_j$ with $j\notin S$ form an increasing subsequence. 
This, and the definition of excedence imply
that $p_j$ is a {\em right-to-left minimum} (that is, smaller than all entries
on its right) if and only if $j\notin {\rm Exc}(p)=S$. 

  Now let $p'= p_n p_{n-1} \cdots p_1$ be the reverse of $p$. 
  Then $p'$ is a 123-avoiding permutation
 having a left-to-right minimum at position $i\leq n$ exactly if
 $n+1-i \notin S$.

There is exactly one 132-avoiding permutation $p''$ which has this
same set of left-to-right minima at these same positions \cite{schmidt}. 
Namely,  $p''$ is obtained by keeping the left-to-right minima of $p'$ fixed, 
and successively placing in the remaining positions,  from left to
right, the smallest available element which does not alter the
 left-to-right minima. 
 We set $s(p)=p''$. Observation 1 then tells us that $i\in {\rm Des} (p'')$
 if and only
 if $n-i\notin S$, in other words, when $i\in \alpha(S)$,
 and so $p''$ belongs indeed to $D^{132}_n(\alpha(S))$.

It is easy to see that $s$ is invertible. Clearly, $p'$
can be recovered from $p''$ as the only 123-avoiding permutation with the
same values and positions of its left-to-right minima as $p''$. (All
entries which are not left-to-right minima are to be written in decreasing 
order). Then $p$ can be recovered as the reverse of $p'$.

The bijections $s \colon E^{321}_n(S) \to D^{132}_n(\alpha(S))$ 
for all the choices of $S \subseteq [n-1]$ produce an order-reversing
bijection from $Q^A_n$ to $P^A_n$.
But $P^A_n$ is self-dual, so the proof is complete.
\end{proof}

\begin{example}\label{eg-sbij}
{\em Take $p=2\ 4\ 1\ 6\ 3\ 5\ 9\ 10\ 7\ 8 \in
E^{321}_{10}(S)$ for $S= \{ 1, 2, 4, 7, 8 \}$.  Then its reversal  
$p'=8\ 7\ 10\ 9\ 5\ 3\ 6\ 1\ 4\ 2$ has left-to-right minima 
8, 7, 5, 3, 1 in positions 1, 2, 5, 6, 8.  We obtain 
$s(p) =  p''= 8\ 7\ 9\ 10\ 5\ 3\ 4\ 1\ 2\ 6$, a permutation in 
$D^{132}_{10}( \{ 1, 4, 5, 7 \})$.
}
\end{example}

\section{The type-B case}\label{sec-B} 

\subsection{The type-B noncrossing partitions}\label{subsec-NCB}

The hyperplane arrangement of the root system of type $B_n$ consists of the
hyperplanes with equations $x_i = \pm x_j$ for $1 \leq i < j \leq n$ and the
coordinate hyperplanes $x_i = 0$, for $1 \leq i \leq n$. 
The subspaces of ${\bf R}^n$ arising as intersections of hyperplanes from
among these can be encoded by partitions of $\{ 1, 2, \dots, n , 
{\overline{1}}, {\overline{2}}, \dots, {\overline{n}} \}$ satisfying the
following properties:  i) if $\{ a_1, \dots, a_k \}$ is a block, then 
$\{ {\overline{a_1}}, \dots , {\overline{a_k}} \}$ is also a block, where
the bar operation is an involution;  and ii) there is at most one block,
called {\em the zero-block}, 
which is invariant under the bar operation.  The collection of such
partitions are the {\em type-B partitions of $[n]$}.  If $1, 2, \dots, n,
{\overline{1}}, {\overline{2}}, \dots, {\overline{n}}$ are placed around a
circle, clockwise in this order, and if cyclically successive elements of the 
same block are joined by chords drawn inside the circle, then, following
\cite{Reiner-NCB}, 
the class of {\em type-B noncrossing partitions}, denoted $NC^B_n$, is the
class of type-B partitions of $[n]$ which admit a circular diagram with no
crossing chords.  Alternatively, a type-B partition is noncrossing if there
are no four elements $a, b, c, d$ in clockwise order around the circle, so
that $a,c$ lie in one block and $b,d$ lie in another block of the partition.
The total number of type-B noncrossing partitions of $[n]$ is 
${ {2n} \choose n}$ (see \cite{Reiner-NCB}).
As in the case of type A, the refinement order on type-B partitions yields a
geometric lattice (in fact, isomorphic to a Dowling lattice with an order-2
group), and the noncrossing partitions constitute a sub-meet-semilattice as
well as a lattice in its own right.  As a poset under the refinement order,
$NC^B_n$ is ranked, with 
${\rm rk}(\pi) =  n - \# ($ of pairs of  non-zero blocks $)$.
For example, $\pi = \{ 1, {\overline{3}}, {\overline{5}} \}, 
\{ {\overline{1}}, 3, 5 \},  \{ 4 \}, \{ {\overline{4}} \}, 
\{ 2, {\overline{2}} \}$ is an element of $NC^B_5$ having 2 pairs of
non-zero blocks and its rank is equal to 3. 
The rank-sizes in $NC^B_n$ are given by $\left( {n \choose k}^2
\right)_{0 \leq k \leq n}$ (see \cite{Reiner-NCB}).  

The numerous properties of $NC^A_n$
which also hold for $NC^B_n$ (as shown in \cite{Reiner-NCB},
\cite{Hersh-NC}), establish the latter as a natural B-analogue. 
In particular, $NC^B_n$ is a self-dual, rank-unimodal, strongly Sperner poset, 
analogously to the properties of $NC^A_n$ of concern in Section \ref{sec-A}.
We now turn to a type-B counterpart of the restricted permutations
considered in the preceding section.

\subsection{A class of pattern-avoiding signed
permutations}\label{subsec-Brestr}

We will view the elements of the hyperoctahedral group $B_n$ as signed
permutations written as words of
the form $b = b_1 b_2 \dots b_n$ in which each of the symbols $1, 2, \dots,
n$ appears, and may or may not be barred.  Thus, the cardinality of $B_n$ is
$n! 2^n$.  
To find a B-analogue of the poset $P^A_n$, we need a subset of $B_n$ whose
cardinality is $\# NC^B_n = { {2n} \choose n }$,  which is characterized via
pattern-avoidance, and over which the distribution of the descent statistic
agrees with the distribution across ranks of the type-B noncrossing partitions 
of $[n]$.  Such a class of signed permutations is 
$B_n(12, {\overline{2}}\ {\overline{1}})$ which appears in
\cite{Si-Bstats}.  We include its description for the reader's convenience.

Consider the elements of $B_n$ which avoid simultaneously the patterns $21$
and ${\overline{2}}\ {\overline{1}}$.  That is, the set of elements $b = b_1
b_2 \cdots b_n \in B_n$ such that there are no indices 
$1 \leq i < j \leq n$ for which i) either both $b_i, b_j$ are barred, or
neither is barred, and ii) $|b_i| > |b_j|$ (the absolute value of a symbol 
means $|a| = a$ if $a$ is not barred, and $|a| = {\overline{a}}$ if $a$ is
barred;  effectively, the absolute value removes the bar from a barred
symbol).
  The following is immediate:  a $(21, {\overline{2}}\ 
{\overline{1}})$-avoiding permutation in $B_n$ is a shuffle of an
increasingly ordered subset $L$ of $[n]$ whose elements we then bar, with its
increasingly ordered complement in $[n]$.  For example, 
$b = {\overline{ 2}} 1 3 {\overline{5}} 4 {\overline{ 6}} 7$ 
is one of 
${7 \choose 3}$ elements of $B_7(21, {\overline{2}}\ {\overline{1}})$ 
associated with the subset $L = \{ 2, 5, 6 \} \subseteq [7]$.
Obviously, summing over the choices of $L$ of cardinality ranging from zero
to $n$ and over the shuffles, it follows that 
\begin{equation}\label{eq-Brestr-card}
\# B_n(21, {\overline{2}}\ {\overline{1}}) = \sum_{k=0}^n
{ {n \choose k}^2 } = { {2n} \choose n} = \# NC^B_n,
\end{equation}
as desired.

Furthermore, the distribution of descents over 
$B_n(21, {\overline{2}}\ {\overline{1}})$ is as desired. 
We say that $b = b_1 b_2 \dots b_n \in B_n$ 
has a {\em descent} at $i$, for $1 \leq i \leq n-1$, if $b_i > b_{i+1}$ with
respect to the total ordering 
$1 < 2 < \cdots < n < {\overline{n}} < \cdots < {\overline{2}} < 
{\overline{1}}$, and that it has a descent at $n$ if $b_n$ is barred.
As usual, the {\em descent set} of $b$, denoted ${\rm Des}(b)$, is the set
of all $i \in [n]$ such that $b$ has a descent at $i$.
For example, for  
$b = {\overline{ 2}} 1 {\overline{3}} {\overline{5}} 4 7 {\overline{ 6}}$ 
we have 
${\rm Des}(b) = \{ 1, 3, 4, 7 \}$.
It is then transparent that if 
$b \in B_n(21, {\overline{2}}\ {\overline{1}})$, then its descent set is
precisely the set of positions occupied by barred symbols. 
In conclusion,

\begin{observation}\label{obs-bij-B2121}
For an element $b$ of the hyperoctahedral group $B_n$, let 
$L(b)$ denote the set of symbols which are barred in $b$, and 
${\rm Des}(b)$ denote the descent set of $b$. 
Then the map $b \mapsto (L(b), {\rm Des}(b))$  gives a bijection between 
the class of restricted signed permutations $B_n(21, {\overline{2}}\
{\overline{1}})$ and ordered pairs of subsets of $[n]$ of equal cardinality.
\end{observation}

\subsection{The poset $P^B_n$}\label{subsec-PBn}

As the B-analogue of the poset of 132-avoiding permutations $P^A_n$ of the
preceding section, we consider the poset $P^B_n$ consisting of 
the $(21, {\overline{2}}\ {\overline{1}})$-avoiding elements of the
hyperoctahedral group $B_n$, with the order relation given by 
$b < b'$ if and only if ${\rm Des}(b) \subset {\rm Des}(b')$.

Based on the preceding discussion and an encoding of type-B noncrossing
partitions appearing in \cite{Reiner-NCB}, one readily obtains the
properties of $P^B_n$ which parallel those of $P^A_n$.

\begin{theorem}
The poset $P^B_n$ of $(12,{\overline{2}}\ {\overline{1}})$-avoiding elements
of the hyperoctahedral group $B_n$, ordered by containment of the descent
set, is an extension of the refinement order on the type-B noncrossing 
partition lattice $NC^B_n$.  The poset $P^B_n$  has the same 
rank-generating-function as $NC^B_n$, therefore it is rank-symmetric and
rank-unimodal, and it is a self-dual and strongly  Sperner poset.
\end{theorem}

\begin{proof}
It is immediate from its definition and Observation \ref{obs-bij-B2121} that 
$P^B_n$ is a ranked poset (namely, ${\rm rk}(b) = \# {\rm Des}(b)$)
and has rank-sizes given by 
$\left( {n \choose k}^2 \right)_{0 \leq k \leq n}$, equal to the rank-sizes
in  $NC^B_n$.
Also, $P^B_n$ is a self-dual poset:  
clearly, if $b'$ is the $(21, {\overline{2}}\ {\overline{1}})$-avoiding
signed permutation which corresponds to the pair 
$( [n] - L(b), [n] - {\rm Des}(b) )$, then the mapping $b \leftrightarrow b'$ 
is an order-reversing involution on $P^B_n$.

Toward checking that 
there is an order-preserving bijection from  $NC^B_n$ to $P^B_n$, we 
first recall a fact from \cite{Reiner-NCB}:
 every partition $\pi \in NC^B_n$ 
can be encoded by a pair $(L(\pi), R(\pi))$ of subsets of $[n]$ whose
cardinality is the number of pairs of non-zero blocks of $\pi$. 
Informally, these sets consist of the Left and Right delimiters of non-zero
blocks when the elements are read in clockwise order (in the circular
diagram of $\pi$).  More precisely, 
if $n=0$ or if $\pi$ has only a zero-block, we set $L=R=\emptyset$.
Otherwise,   
$\pi \in NC^B_n$ has some non-zero block consisting of cyclically
consecutive elements in its diagram.  If such a block consists of 
$j_1, j_2, \dots, j_k$ in clockwise order, then $|j_1|$ belongs to $L(\pi)$
and $|j_k|$ belongs to $R(\pi)$.  By deleting this block and
its image under barring, a type-B noncrossing partition of
$[n-k]$ is obtained and the construction of the sets $L(\pi)$ and $R(\pi)$ 
is completed by repeating this process as long as non-zero blocks arise. 
For instance, if $\pi =  
 \{ 1, {\overline{3}}, {\overline{5}} \}, 
\{ {\overline{1}}, 3, 5 \},  \{ 4 \}, \{ {\overline{4}} \}, 
\{ 2, {\overline{2}} \}$, then $L(\pi) = \{ 3, 4 \}$ and 
$R(\pi) = \{ 1, 4 \}$.

Now suppose that $\pi < \pi'$ in $NC^B_n$, and that this is a covering relation
(i.e., ${\rm rk}(\pi') = {\rm rk}(\pi) + 1$).
Then there exist $l \in L(\pi)$ and $r \in R(\pi)$ such that 
$L(\pi') = L(\pi) - \{ l \}$ and $R(\pi') = R(\pi) - \{ r \}$, 
as a result of the merging of blocks entailed by the covering relation.
Thus it is clear that if $\pi \in NC^B_n$ is mapped to the signed
permutation $b \in P^B_n$ with the property that 
$(L(b), {\rm Des}(b)) = (L(\pi), R(\pi))$, then one obtains an order-reversing
embedding of $NC^B_n$ into $P^B_n$.   Combining this with the 
self-duality of $P^B_n$ we obtain the desired embedding of $NC^B_n$ into
$P^B_n$.

Finally, the strong Sperner property of $P^B_n$ follows as in type A, from the 
strong Sperner property of $NC^B_n$ (see \cite{Reiner-NCB}) and the
 rank-preserving
embedding of $NC^B_n$ into $P^B_n$. 
\end{proof}

\subsection{A poset based on type-B excedences}\label{subsec-excB}

As in the type-A case, there is a self-dual poset 
of $\# NC^B_n$  restricted signed permutations ordered by containment of 
the set of
excedences.   In fact, there is more than one definition of the excedence
statistic in the literature, in the case of the hyperoctahedral group. 
We briefly mention two possibilities considered in \cite{Steing}. 

Given a signed permutation $b$, let $k$ be the number of symbols which are
{\em not} barred in $b$.  
We associate to $b$ an $(n+1)$-permutation $\sigma(b)$ by setting
$\sigma(b)_{n+1} = k+1$ and, for $1 \leq i \leq n$, letting  
$\sigma(b)_i = j$ if $b_i$ is the $j$th smallest among the symbols 
$b_1, b_2, \dots, b_n, n+1$ with respect to the linear ordering 
$1 < 2 < \cdots < n < n+1 < {\overline{1}} < {\overline{2}} < \cdots <
{\overline{n}}$. 
For example, if $b = 1\ {\overline{3}}\ 2\ 4\ 5\ {\overline{6}}\ 
{\overline{8}}\ 7$, then $\sigma(b) = 1\ 7\ 2\ 3\ 4\ 8\ 9\ 5\ 6$.
Now, the excedence set of $b$ is defined to be that of $\sigma(b)$.
It turns out \cite{Gal-pers} that for 
$b \in B_n(21, {\overline{2}} \ {\overline{1}})$ this
definition makes the excedence set coincide with the descent set for each
$b$.  Therefore, this leads to the poset $P^B_n$ again.

An alternative definition for excedences of ``indexed permutations'' appears
in \cite{Steing}.  Specialized to the hyperoctahedral group it is the 
following.

\begin{definition}\label{def-Bexc2}
If $b \in B_n$, its {\em excedence set} is the union of 
the sets $S(b)$ and $F(b)$, where
$S(b)$ is the set of excedences computed in the symmetric group for the
permutation $|b_1| |b_2| \dots |b_n|$ 
obtained by removing all bars from the symbols in $b$,
and $F(b) = \{ i \ \colon \ b_i = {\overline{i}} \}$, 
the set of barred fixed points of $b$.  
\end{definition}
Thus, for  
$b = 1\ {\overline{3}}\ 2\ 4\ 5\ {\overline{6}}\ {\overline{8}}\ 7$
we obtain excedences at 
$\{ 2, 6, 7 \}$ by either of the two definitions. 
But  $b = {\overline{1}} 3 {\overline{2}}$ has excedences at $\{ 1, 3 \} $
by the first definition (based on $\sigma(b) = 3\ 1\ 4\ 2$),
and $\{ 1, 2 \}$ if the second definition is adopted. 

For the remainder of this section, we work with the notion of
excedence as in Definition \ref{def-Bexc2}.  

\begin{proposition}\label{prop-QBn-sdual}
Let $Q^B_n$ denote the poset of 
$(21, {\overline{2}}\ {\overline{1}})$-avoiding signed permutations in 
$B_n$, ordered by containment of their excedence set. 
The poset $Q^B_n$ is self-dual.
\end{proposition} 

\begin{proof} Let $b\in B_n$ and $b'$ be the reverse of $b$.  Let
$b''$ be the ``barred complement'' of $b'$, that is,  
$|b''_i| = n+1- |b'_i|$, and  $b''_i$ is barred if and only if 
$b'_i$ is not barred.

 Then it is straightforward to verify that $i\in S(p'')\cup F(p'')$
if and only if  $i\notin S(p)\cup F(p)$. Therefore, the reverse complement
operation reverses the inclusion of excedence sets for signed permutations.
(Thus, the entire hyperoctahedral group $B_n$ ordered by containment of the
excedence set is a self-dual poset.)
But, clearly, this involution preserves the $(21, {\overline{2}}\
{\overline{1}})$-avoidance property, and thus $Q^B_n$ is self-dual.
\end{proof}

By \cite{Steing}, the rank generating function of $Q^B_n$ is equal to that of 
$P^B_n$.  Therefore it is natural to ask whether the posets $P^B_n$   
and $Q^B_n$ are isomorphic, just as their type-A counterparts are
(Proposition \ref{prop-isoA}).  The answer in this case is negative. 
Indeed, if $n=3$ it is straightforward to verify that all atoms of $P^B_3$ 
are covered by six elements, while the atom 
${\overline{1}}\ 2\ 3$ of $Q^B_3$ is covered by seven elements 
(namely, ${\overline{1}}\ {\overline{2}}\ 3, 
{\overline{1}} \ 2 \ {\overline{3}}, 
{\overline{1}} \ 3 \ {\overline{2}}, 
2 \ 3 \ {\overline{1}},
{\overline{2}} \ {\overline{3}} \ 1,
2 \ {\overline{1}} \ {\overline{3}}$, and  
${\overline{2}} \ 1 \ {\overline{3}}$).

\section{Remarks and questions for further investigation}\label{sec-Further}

\begin{enumerate}

\item{ {\em Is $NC^B_n$ a subposet of $Q^B_n$?} 
We do not know whether the lattice of type-B noncrossing partitions 
can be embedded in the poset $Q^B_n$ of 
$(21, {\overline{2}}\ {\overline{1}})$-avoiding signed permutations ordered
by their excedence set of definition \ref{def-Bexc2}.
}

\item{  {\em Self-duality of $NC^A_n$ and $NC^B_n$
extending to self-duality for $P^A_n$ and $P^B_n$.}

We have seen that {\em each} of the posets $NC^A_n$ and $P^A_n$ is self-dual 
and  that $NC^A_n$ is a subposet of $P^A_n$. The same is true for the pair 
$NC^B_n$, $P^B_n$.  Both for type A and for type B one can exhibit an
order-reversing involution on the larger poset which restricts to an
order-reversing involution on the smaller one. 

We first construct such  an involution $g$ for $NC^A_n$ 
which will be similar, though not identical,
 to the involution defined in \cite{SiU}.

Write the elements $1,2,\cdots ,n$ clockwise around a circle, and write
elements $1',2',\cdots ,n'$ interlaced in counterclockwise order,
so that $1'$ is between $1$ and  $n$, $2'$ is between $n$ and $n-1$, and so
on, $i'$ is between $n+2-i$ and $n+1-i$. For $\pi \in NC^A_n$, join by
chords -- as usual -- cyclically successive (unprimed) elements
belonging to the same block of $\pi$. Then define $g(\pi)$ to be
the coarsest noncrossing partition on the elements $1',2',\cdots ,n'$ so
that the chords joining primed elements of the same block do 
not intersect the chords of $\pi$. See Figure 2 for an example.

\begin{figure}[ht]
 \begin{center}
  \epsfig{file=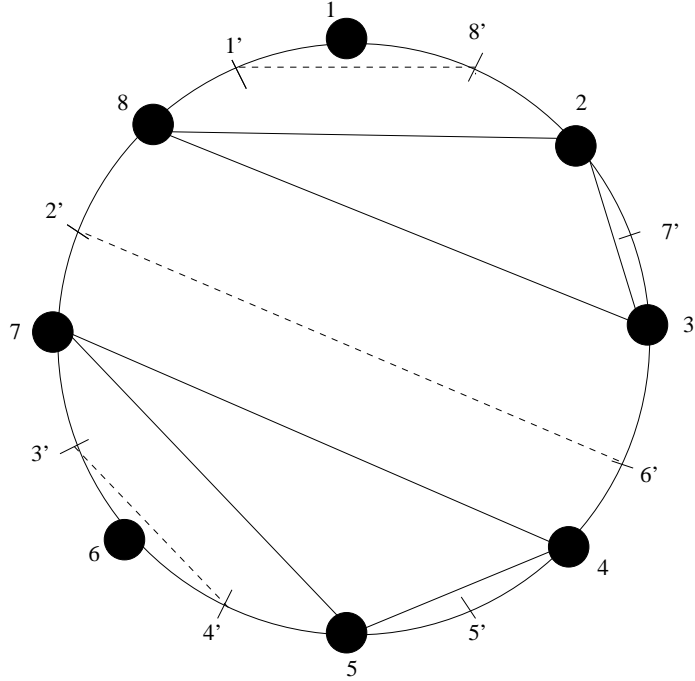}
  \caption{The partition $\pi=(\{1\},\{2,3,8\},\{4,5,7\},\{6\})$ and its image
  $g(\pi)$}
\label{duality}
 \end{center}
\end{figure}

The map $g$ is certainly a bijection, and it is
 order-reversing in $NC^A_n$ since merging two blocks of $\pi$ subdivides a
 block of $g(\pi)$. We claim that $g$ is also order-reversing on $P^A_n$.
 To see this, 
observe that for any $i>1$, the element $i$ is the smallest in its block 
in $\pi$ if and
only if the element $(n+2-i)'$ is {\em not} the smallest in its block in 
$g(\pi)$.
Indeed, the definition of $g$ implies that exactly one of
$i$ and $(n+2-i)'$ can be connected to smaller elements by a chord.
Therefore, $g$ takes the set of block-minima (not equal to 1) of $\pi$
into its reverse complement in $[n+2]$, so $g$ is indeed order-reversing on
 $P^A_n$.

In the type-$B$ case, one can obtain an analogous bijection $h$
 in a similar way: take the circular
(clockwise) representation of $\pi$, then write the elements $1',2',\cdots,
n',\bar{1}',\bar{2}',\cdots ,\bar{n}'$ so that the primed numbers interlace
the unprimed, placing $1'$  between $1$ and $\bar{n}$ and continuing 
counterclockwise.  For $\pi \in NC^B_n$, define $h(\pi)$ as
above, that is, as the unique coarsest partition on the primed set
whose chords do not cross those of $\pi$. 
Then $h$ is certainly an order-reversing bijection of
$NC^B_n$, and as above, it reverses the containment of the sets $L(\pi)$
and $R(\pi)$, so it does extend to an order-reversing bijection of $P^B_n$.
}

\item{ {\em The M\"obius function and order complexes of $P^A_n$ and
$P^B_n$.}

It is easy to write an expression for the number 
$c_m(P^B_n )$ of chains ${\hat 0} < b^1 < b^2 < \cdots < b^m < {\hat 1}$ 
of length $m+1$ in $P^B_n$, for $m \ge 0$.
Of course, $c_0(P^B_n) = 1$, and 
\begin{equation}\label{eq-PBchains}
c_m(P^B_n) = \sum_{0 <k_1 <k_2 < \cdots < k_m < n}
   { { n \choose {k_1}} {n \choose {k_2}} \cdots {n \choose {k_m}}
     { {n!} \over { k_1! (k_2 - k_1)! \cdots (k_m - k_{m-1})! (n - k_m)! }}}
\end{equation}
since under the correspondence 
$b^i \leftrightarrow (L(b^i), {\rm Des}(b^i))$
a chain in $P^B_n$ corresponds to an $m$-tuple of subsets $(L(b^i))$ and a
chain of subsets ${\rm Des}(b^1) \subset {\rm Des}(b^2) \subset \cdots 
\subset {\rm Des}(b^m)$ of $[n]$, with $\# L(b^i) = \# {\rm Des}(b^i) =
k_i$.
In turn, this leads to an expression for the M\"obius function of $P^B_n$, 
$\mu_{P^B_n}({\hat 0}, {\hat 1}) = \sum_{m \ge 0}{ (-1)^{m-1} c_m(P^B_n)}$.

These expressions can be regarded as partial success with the computation 
of the zeta polynomial and the M\"obius function.  It would be
interesting to elucidate further the question of these invariants for 
$P^A_n$ and $P^B_n$,  and to describe the order complexes of these posets. 
}

\item{ {\em Other posets of combinatorial objects with similar properties.} 

The behavior of noncrossing partitions and restricted permutations suggests
the following question:  what other combinatorial objects 
admit a natural partial order which is self-dual and possibly,
has other nice properties?  A natural candidate is the class of {\em two-stack
sortable permutations} \cite{doron}. It is known \cite{schaeffer} 
that there are as many of them with $k$ descents as with
$n-1-k$ descents.  However, the poset obtained by the descent ordering is
not self-dual, even for $n=4$, so another ordering is needed. 

Similarly, the type-D noncrossing partitions and the interpolating
BD-noncrossing partitions do not, in general, form self-dual posets when
ordered by refinement (see \cite{Reiner-NCB}).   However, it may be
interesting to find corresponding classes of pattern-avoiding elements in
the Weyl group for type D, along with an order-preserving embedding 
$NC^D_n \to P^D_n$ analogous to the type-A and B cases.}

\end{enumerate}

\end{document}